# Constrained and Preconditioned Stochastic Gradient Method

Hong Jiang, Gang Huang, Paul Wilford and Liangkai Yu

*Abstract*—We consider stochastic approximations which arise from such applications as data communications and image processing. We demonstrate why constraints are needed in a stochastic approximation and how a constrained approximation can be incorporated into a preconditioning technique to derive the preconditioned stochastic gradient method (PSGM). We perform convergence analysis to show that the PSGM converges to the theoretical best approximation under some simple assumptions on the preconditioner and on the independence of samples drawn from a stochastic process. Simulation results are presented to demonstrate the effectiveness of the constrained and preconditioned stochastic gradient method.

## I. INTRODUCTION

STOCHASTIC approximations have many applications in data communications and image processing, including channel equalization, digital predistortion for high power amplifiers and approximation of camera response function [1]-[8]. They are also used in the fields of medical devices [9], machine learning and data miming [10]. The stochastic approximation of interest in this paper can be characterized by 1) the function to be approximated is not known; 2) the input and output of the function can be observed with samples, and the input samples are drawn from a stochastic process of a certain probability density function, and its distribution can be observed, but cannot be controlled or altered; and 3) there is an unlimited supply of input and output samples for observations, but there may be a limit on how many samples one can observe at one time.

Stochastic approximations usually involve approximating the unknown function by a linear combination of some basis functions. The coefficients of the basis functions in the linear combination satisfy a system of linear equations which are determined by the unknown function and the distributions of the input and output signals.

Stochastic gradient descent method has been widely used for stochastic approximation [11]-[13]. In the stochastic gradient method (SGM), the approximation is computed in an iterative process in which each step is marched in the gradient descent direction of a cost function evaluated at a set of samples.

Two issues arise frequently in the process of SGM. First, a solution computed at a step with a particular set of samples may have an undesirable behavior, e.g., with very large derivative, or the linear system can even be ill-posed with a singular matrix. Therefore, it may be preferred or necessary to impose additional constraints on the solution of SGM so that it has a desirable behavior, e.g., smaller derivatives, or the linear system becomes well-posed with an invertible matrix. Secondly, it is well-known that gradient descent method may converge slowly, and some techniques are needed to speed up the convergence.

### A. Main contributions of this paper

In this paper, we exemplify the need for constrained stochastic approximations and propose a few effective constraints which are simple to implement. We also use a preconditioning technique [14] to speed up the convergence of SGM and propose the preconditioned stochastic gradient method (PSGM) for stochastic approximations.

Furthermore, we establish a link between the constrained stochastic approximation and the preconditioning, and show that the constraints can be incorporated as part of the preconditioner in PSGM.

Finally, we provide a theoretical analysis to show that the PSGM converges under some simple assumptions on the preconditioner and the stochastic process. The convergence can be proved with the simple assumption that the sample sets are drawn identically and independent from one set to another, without any assumption about the distribution of the random samples.

### B. Related work

Stochastic approximations have been analyzed extensively in [1]. Convergence analysis of stochastic conjugate gradient method were performed in [2][5], but the stochastic conjugate gradient method has more complexity than the PSGM of this paper. In addition, the analysis of [2] also required more restrictive assumptions on the distribution of certain random variables derived from the random samples. Analysis of SGM for linear adaptive filters was given in [15] and [16], but the results there were obtained for some known distributions of the stochastic process. Preconditioning has also been used in combination with stochastic gradient descent in other fields, see for example, [9][10].

### C. Organization of the paper

In section II, the stochastic approximation problem of interest of this paper is posed. This is followed by a few examples given in section III to demonstrate the need for constrained approximations and to establish the link between constrained approximation and the preconditioning. The preconditioned stochastic gradient method is formally introduced in section IV, and it is analyzed in section V. Simulation results are reported in section VI.

## II. PROBLEM FORMULATION

### A. Theoretical best approximation

Let $X$ be a random variable or a stationary stochastic process, $X = X(t)$, with the sample space an interval of $\Re$ and the probability density function $\rho(x)$. Let $f(x)$ be an un-The authors are with Bell Labs, Alcatel-Lucent, Murray Hill, NJ, USA.

known deterministic function defined on the sample space of $X$. Function $f(x)$, although unknown, may be observed through the stochastic process

$$Y = f(X). \quad (1)$$

The objective is to find an approximation of $f(x)$ by using a linear combination of some basis functions, shown in Fig 1(a).

Let $\Phi(x)$ be a row vector of length $M$ given by

$$\Phi(x) = [\phi_1(x), \phi_2(x), ..., \phi_M(x)], \quad (2)$$

where $\phi_j(x), j = 1,...,M$, are linearly independent basis functions [17] defined in the sample space of $X$. Any function in the subspace spanned by $\Phi(x)$ can be written as

$$u(x) = \Phi(x) \cdot u = \sum_{j=1}^{M} u_j \phi_j(x), u = [u_1, u_2, ..., u_M]^T \in \Re^M. \quad (3)$$

In (3), variable $u$ is used to denote both a function and a column vector, but no confusion should arise because it is clear from the context whether $u$ refers to a function or a vector.

The approximation problem is then formulated as finding a vector $\hat{u}$ such that the function $\Phi(x) \cdot \hat{u}$ best approximates $f(x)$. The best approximation is in the sense as measured by the expected value of the squares of the residual. More precisely, for any function $g(x)$, let $E(\cdot) = E_\rho(\cdot)$ be the expectation function defined as

$$E(g(X)) = \int g(x) \rho(x) dx. \quad (4)$$

Then, the best approximation $\hat{u}$ is defined by

$$\hat{u} = \arg\min_{u \in \Re^M} E((\Phi(X) \cdot u - f(X))^2). \quad (5)$$

From (5), it can be easily shown that the best approximation $\hat{u}$ is the solution to the linear equation

$$Au = b, \quad (6)$$

where matrix $A$ and vector $b$ are given by

$$\begin{aligned} b &= [b_1, ..., b_M]^T, b_i = E(\phi_i f), i = 1, ..., M, \\ A &= [a_{ij}]_{M \times M}, a_{ij} = E(\phi_i \phi_j), i, j = 1, ..., M. \end{aligned} \quad (7)$$

The solution $\hat{u}$ to (6) provides the theoretical best approximation function $\hat{u}(x)$ which is given by

$$\hat{u}(x) = \Phi(x) \cdot \hat{u} = \sum_{j=1}^{M} \hat{u}_j \phi_j(x). \quad (8)$$

Finding the solution to (6) would otherwise be trivial if it were not for the fact that the right hand side $b$, give in (7), cannot be computed because $f(x)$ is not known explicitly. In practice, (6) is approximated by a sequence of equations which are computed by using samples from the stochastic process $X$, hence giving rise to the stochastic approximation.

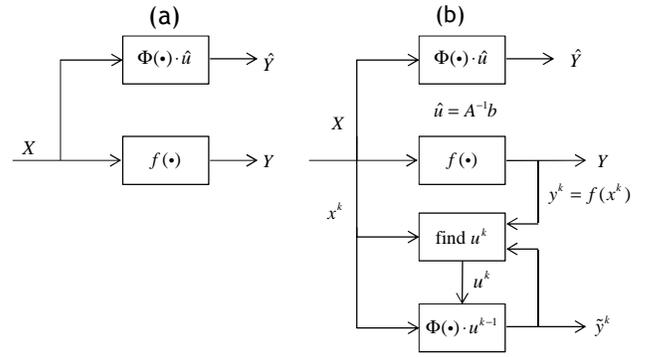

Fig 1 (a) Approximation of $f(\cdot)$. (b) Stochastic approximation

### B. Stochastic approximation

Let $x^k, k = 1, 2, ...$, be a sequence of sets of samples drawn from $X$ as

$$x^k = [x_1^k, ..., x_N^k]^T \in \Re^N, k = 1, 2, ..., \quad (9)$$

That is, for each $k = 1, 2, ...$, there are $N$ samples drawn from $X$, and they are given by $x_n^k, n = 1, 2, ..., N$. Then the values of the function $f(x)$ can be observed at the samples, and assuming no noise, they are given by

$$y_n^k = f(x_n^k), \ n = 1, ..., N, k = 1, 2, ... \quad (10)$$

The stochastic approximation is an iterative process in which a new approximation $u^k$ is computed from previous $u^{k-1}$ and the new observation $y_n^k$ of (10). The new coefficient $u^k$ is updated from $u^{k-1}$ so that the function $\Phi(x) \cdot u^k$ best matches $f(x)$ when evaluated at the sample set $x^k$:

$$\Phi(x_n^k) \cdot u^k = f(x_n^k) = y_n^k, \ n = 1, ..., N. \quad (11)$$

The new coefficient can be expressed as

$$u^k = u^{k-1} + \Delta u^k. \quad (12)$$

If we define matrix $\Phi_k$ as $\Phi(x)$ evaluated at the samples, i.e.,

$$\Phi_k = \left[\phi_j(x_i^k)\right]_{i=1,...,N, j=1,...,M} \in \Re^{N \times M}, \quad (13)$$

and substitute (12) into (11), the update $\Delta u^k$ is given by

$$\Phi_k \Delta u^k = r^k, \text{ where } r^k = y^k - \Phi_k u^{k-1}. \quad (14)$$

Equations (12) and (14) define the stochastic approximation as an iterative process. At step $k$ of the iteration, a set of $N$ samples, $x^k$, are drawn to compute residual $r^k$ in (14). Then the update, $\Delta u^k$, is determined by requiring it to satisfy (14). The new coefficient $u^k$ is computed from (12), and the iteration then proceeds to the next step. This process is illustrated in Fig 1(b).

### C. Least squares solution

One implementation of the stochastic approximation (12) is to use the least squares solution in (14), which is an overdetermined system if the number of samples is larger than the number of basis functions, i.e., if $N > M$. The least squares solution to (14) is given by

$$A_k \Delta u^k = b^k - A_k u^{k-1}, \quad (15)$$

where $A_k$ and $b^k$ are given by

$$A_k = \Phi_k^T \Phi_k = [a_{ij}^k]_{M \times M}, \ a_{ij}^k = (1/N)\sum_{n=1}^{N} \phi_i(x_n^k)\phi_j(x_n^k), \quad (16)$$

$$b^k = \Phi_k^T y^k = [b_1^k,...,b_M^k]^T, \ b_i^k = (1/N)\sum_{n=1}^{N} \phi_i(x_n^k)y_n^k. \quad (17)$$

Note that for any set of samples $x^k$, (16) and (17) can be considered as approximations of $A$ and $b$ in (7). It is always assumed that the sets of samples $x^k$ are chosen properly so that matrix $A_k$ is symmetric and positive definite, i.e.,

$$A_k^T = A_k, \ x^T A_k x > 0, \ \text{for all } x \neq 0. \quad (18)$$

Therefore, the update $\Delta u^k$ can be found from (15), and (12) is now equivalent to

$$u^k = A_k^{-1} b^k. \quad (19)$$

### D. Challenges

The following are a few issues with using (19) in the stochastic approximation process of Fig 1(b).

1. Matrix $A_k$, although symmetric and positive definite, may be ill-conditioned. That is, the condition number, which is defined as the ratio of the largest eigenvalue over the smallest eigenvalue, is very large

$$\kappa(A_k) \triangleq \lambda_{max} / \lambda_{min} \gg 1. \quad (20)$$

A large condition makes solution of (15) or (19) sensitive to noise and round-off errors. If iterative method is used, the convergence may be slow due to large condition number.

2. Solution $u^k$ of (19) is only the best approximation over $N$ samples $x_n^k, n = 1,...,N$. At the points other than these samples, the function $u^k(x) = \Phi(x) \cdot u^k$ may behave badly. For example, at a point $z > \max_n\{x_k^n\}$, the value $\Phi(z) \cdot u^k$ may be unfavorably large.

3. The dimension of $A_k$ may be large, i.e., $M \gg 1$, in which case, computing $A_k^{-1} b_k$ may be too costly in both the number of operations and memory requirement.

4. Solutions $u^k$ of (19) do not converge to the best approximation $\hat{u}$. In fact, if the sets of samples $x^k$ are drawn randomly and independently for each $k$, then all solutions are statistically the same, and no matter how many sample sets are used, the solution will not get any better.

5. Even if solutions $u^k$ of (19) are averaged, it is still not clear whether the average converges to $\hat{u}$, because the limit $\lim_{K \to \infty}(1/K)\sum_{k=1}^{K} A_k^{-1}b^k$ may not equal $A^{-1}b$, or even exist, under reasonable assumptions.

### III. CONSTRAINED STOCHASTIC APPROXIMATION

In many applications, the approximation $u^k$ computed from a particular set of samples may exhibit undesirable behaviors and it is necessary to impose additional constraints on the approximation in order to alleviate the problem. For example, it may be desirable to constrain the solution so that the change in derivatives from one sample set to anther is small (see below). In this section, we discuss how the constraints can be introduced into the stochastic approximation process with a few examples. This also serves as a motivation for preconditioned method to be introduced later.

### A. Polynomial approximation

In polynomial approximations, the basis functions $\Phi(x)$ are made of polynomials, or orthogonal polynomials [18], so that $\phi_j(x)$ is a polynomial of degree $j$. For each set of samples $x^k, k = 1,...$, we want to find coefficients $u^k$ so that the polynomial $\Phi(x) \cdot u^k$ matches the observed output of the function $y^k$ when evaluated at the samples. By following the stochastic process of Fig 1(b) as described before, the coefficients may be computed from (12) and (15), or equivalently (19).

As is well known, the coefficients thus computed results in a polynomial $\Phi(x) \cdot u^k$ that may behave badly at the points other than the samples $x^k$. For example, as will be shown in the simulation section, if the interval of approximation is [0,1], and if the largest sample is significantly less than 1, then the derivative of $\Phi(x) \cdot u^k$ may be very large in the interval between the largest sample and 1, i.e.,

$$\|\frac{d}{dx}\Phi(z) \cdot u^k\| \gg 1, \ z \in (\max_n\{x_n^k\},1]. \quad (21)$$

In such a case, it is desirable to place a constraint on the approximating polynomial so that it has a small derivative. Such a constraint may be constructed as follows. First, we select a set of samples where the derivative is to be evaluated.

$$z_i \in (\max_n\{x_n^k\},1], i = 1,...,N_c. \quad (22)$$

Then we impose the condition that the change in derivative from one sample set to another be zero at the selected points:

$$\frac{d}{dx}\Phi(z_i) \cdot \Delta u^k = 0, \text{ or}$$
$$D\Delta u^k = 0, \text{ where } D = \frac{d}{dx}\Phi(z_i) = \left[\phi_j'(z_i)\right]_{N_c \times M} \in \Re^{N \times M}. \quad (23)$$

Now to find the constrained approximation, the update $\Delta u^k$ is required not only to match the observed function values as given in (14), but also to satisfy the zero derivative condition of (23). A new equation for $\Delta u^k$ is obtained when constraint (23) is combined with the original equation (14), as follows

$$[\Phi_k \ \ D]^T \Delta u^k = [r^k \ \ 0]^T. \quad (24)$$

Using the least squares solution for (24), the new updated coefficients is given by

$$u^k = u^{k-1} + (A_k + D^T D)^{-1}(b^k - A_k u^{k-1}). \quad (25)$$

More generally, we introduce a parameter $\gamma > 0$ to control the extent of the constraint in (23), and a damping factor $\mu > 0$ to control the rate of the update so that (25) becomes

$$u^k = u^{k-1} + \mu(A_k + \gamma D^T D)^{-1}(b^k - A_k u^{k-1}). \quad (26)$$

Then a question arises: how are coefficients $u^k$ computed from the constrained update (26) related to the best approximation $\hat{u}$ defined in (5)?

### B. Approximation by piece-wise constant

A Look-up-table (LUT) can be used to represent a function. To use LUT, the sample space $\Gamma$ is first discretized into small disjoint intervals as

$$\Gamma = \bigcup_{j=1}^{M} \Gamma_j, \quad \Gamma_j = [\chi_{j-1}, \chi_j], j=1,...,M. \quad (27)$$

Then an LUT is an array $u = [u_1,...,u_M]^T$ which can be used to represent a piecewise constant function $g_u(x)$ defined by

$$g_u(x) = u_j, \text{ if } x \in \Gamma_j. \quad (28)$$

The basis functions for piecewise functions are given by

$$\phi_j(x) = \phi_{\Gamma_j}(x) = \begin{cases} 1, & x \in \Gamma_j \\ 0, & x \notin \Gamma_j \end{cases}, j=1,...,M. \quad (29)$$

When using an LUT, the number of basis functions, $M$, is usually very large, e.g., $M = 1024$ if a 10 bit LUT is used. Similar to the polynomial approximation, for a given set of samples $x^k$, an LUT $u^n$ can be computed from (12) and (15). However, even if the number of samples, $N$, is large, there is no guarantee that the samples fall into each and every subinterval $\Gamma_j$, and therefore, for any sample set $x^k$, there may exist some subintervals $\Gamma_j$ into which none of the samples in $x^k$ falls. Consequently, the matrix $\Phi_k$ in (14) may not have the full rank. Indeed, if none of the samples belongs to $\Gamma_j$, then column $j$ of matrix $\Phi_k$ is zero. That is, coefficient $\Delta u_j^k$ is not involved in the equation (14) or (15), and therefore, (14) or (15) is underdetermined. In order to find a solution, additional constraints are needed. One constraint can be, for example, the requirement that the update to the LUT is smooth, so that the difference in $\Delta u^k$ is small, i.e.,

$$\Delta u_j^k - \Delta u_{j-1}^k \approx 0, j = 2,...,M. \quad (30)$$

Equation (30) suggests the following constraint to be used

$$D\Delta u^k = 0, \text{ where } D = \begin{bmatrix} -1 & 1 & & \\ & \ddots & \ddots & \\ & & -1 & 1 \end{bmatrix}_{(M-1) \times M}. \quad (31)$$

Applying constraint (31) to (14) leads to (24), and the update (26) can be used for its least squares solution to find the LUT $u^k$ which exhibits smoothness.

### C. Functions with memory

Some functions exhibit memory effect. This happens often in communication systems in which the random variable models a time-varying signal through a channel with memory, and the set of samples are the signal sampled at time instances of equal interval. As an example, we consider a function $f(x)$ that has a memory of two time intervals, although the analysis applies to cases with any memory delays. In this example, $f(x_n)$ may be best approximated by a function of the form

$$u(x_n) + v(x_{n-1}) + w(x_{n-2}). \quad (32)$$

Each of the functions $u(x), v(x), w(x)$ is a linear combination of the same basis functions $\Phi(x)$, i.e.,

$$u(x) = \Phi(x) \cdot u, \quad v(x) = \Phi(x) \cdot v, \quad w(x) = \Phi(x) \cdot w, \quad (33)$$

for some vectors $u, v, w$. The basis functions may be either polynomials or piecewise constant functions, but in the following discussion, the piecewise constant functions are used because the resulting linear equation has a large dimension worthy paying more attention.

In the stochastic process of Fig 1(b), for given coefficients $[u^{k-1}, v^{k-1}, w^{k-1}]^T$, the new coefficients are found by

$$[u^k, v^k, w^k]^T = [u^{k-1}, v^{k-1}, w^{k-1}]^T + [\Delta u^k, \Delta v^k, \Delta w^k]^T, \quad (34)$$

where the update $[\Delta u^k, \Delta v^k, \Delta w^k]^T$ is computed so that the new approximating function matches the observed $f(x)$. As before, the update is to satisfy

$$\Phi_k^0 \Delta u^k + \Phi_k^1 \Delta v^k + \Phi_k^2 \Delta w^k = r^k,$$
$$\text{where } \Phi_k^p = [\phi_j(x_{i-p}^k)]_{i=1,...,N, j=1,...,M} \in \Re^{N \times M}, p = 0, 1, 2. \quad (35)$$

A constraint may be imposed on each of $\Delta u^k, \Delta v^k, \Delta w^k$, so equation (35) combined with the constraints becomes

$$\begin{bmatrix} \Phi_k^0 & \Phi_k^1 & \Phi_k^2 \\ \sqrt{\gamma}D & & \\ & \sqrt{\gamma}D & \\ & & \sqrt{\gamma}D \end{bmatrix} \begin{bmatrix} \Delta u^k \\ \Delta v^k \\ \Delta w^k \end{bmatrix} = \begin{bmatrix} r^k \\ 0 \\ 0 \\ 0 \end{bmatrix}. \quad (36)$$

In (36), $\sqrt{\gamma}D$ represents the operator for constraint, such as that given in (31). Note that the parameter $\gamma$ is introduced in (36) to control the extent of the constraint. Substituting the least squares solution of (36) into (34), the new coefficients are given by

$$[u^k \ v^k \ w^k]^T = [u^{k-1} \ v^{k-1} \ w^{k-1}]^T + (A_k + \gamma C)^{-1} [\Phi_k^{0T} r^k \ \Phi_k^{1T} r^k \ \Phi_k^{2T} r^k]^T. \quad (37)$$

In (37), $A_k$ and $C$ are given by

$$A_k = [\Phi_k^{iT} \Phi_k^j]_{i=0,1,2, j=0,1,2}, \quad C = \begin{bmatrix} D^T D & & \\ & D^T D & \\ & & D^T D \end{bmatrix}. \quad (38)$$

As noted before, for piecewise constant basis functions, matrix $A_k$ may be singular if the samples $x^k$ are not drawn from all subintervals $\Gamma_j$. On the other hand, for any $\gamma > 0$, the matrix $A_k + \gamma C$ becomes nonsingular if $D$ is the difference operator defined in (31). That is, adding the constraint makes an ill-posed problem well-posed.

Note that each of $\Phi_k^{pT} \Phi_k^p, p = 0, 1, 2$, is a diagonal matrix of dimension $M \times M$, because $\Phi(x)$ is given by the piecewise constant basis functions of (29). However, $A_k$ of (38) is not a diagonal matrix. Also the dimension of $A_k$ is very large; it is $3M \times 3M$, having more than nine million entries if 10 bit LUTs are used. On the other hand, the matrix $C$ is a tridiagonal matrix if the smoothness constraint of (31) is used.

Therefore, it is not desirable to have $A_k$ involved in (37) because computing the last term in (37) may be costly since $A_k$ is large and non-sparse. Instead of $A_k$, we want to use a simpler matrix $B$, which approximates $A_k$, in (37) to compute the update as

$$[u^k \ v^k \ w^k]^T = [u^{k-1} \ v^{k-1} \ w^{k-1}]^T + (B + \gamma C)^{-1} [\Phi_k^{0T} r^k \ \Phi_k^{1T} r^k \ \Phi_k^{2T} r^k]^T. \quad (39)$$

For example, $B$ can be the diagonal matrix of $A$ given in (7), i.e., $B = diag(A)$. Note that $A$ of (7) is computable because it does not involve the unknown function $f(x)$.

### D. Canonical form for constrained approximation

The above three examples show that in a constrained stochastic approximation process of Fig 1(b), the approximation can be computed by the iterative process given by

$$u^k = u^{k-1} + \mu(B + \gamma C)^{-1}(b^k - A_k u^{k-1}). \qquad (40)$$

In (40), $\mu > 0, \gamma > 0$ are some positive constants, $B$ is a matrix that is an approximation of $A_k$ and $C$ is a symmetric positive semi-definite matrix related to the constraint. Although $B$ depends on $k$ in some previous examples, we require $B$ to be independent of $k$ in the canonical form (40) for the reasons to become clear in the next section.

In the rest of the paper, we will analyze how $u^k$ is related to the theoretical best approximation $\hat{u}$ of (5), and in particular, under what condition $u^k$ converges to $\hat{u}$.

We note that in (40), the matrix $B$ is introduced to replace $A_k$ in computing $\Delta u^k$, in order to reduce the complexity involved in operating on $A_k$. The motivation for matrix $\gamma C = \sqrt{\gamma} D^T D$ is to impose a constraint

$$\sqrt{\gamma} D \Delta u^k = 0, \qquad (41)$$

on the update $\Delta u^k$. Although the constraint (41) does affect the intermediate solution $u^k$ by, for example, introducing smoothness in $\Delta u^k$, it does not affect the final limit, because as we will show, $u^k$ converges to $\hat{u}$ regardless of the constraint. Because the constraint (41) does not affect the final limit, it is called soft constraint.

## IV. PRECONDITIONED STOCHASTIC GRADIENT METHOD

### A. Stochastic gradient method

The coefficient given in (19) is the solution which minimizes the following quadratic form

$$H_k(u) = (\Phi_k u - y^k)(\Phi_k u - y^k)^T, \qquad (42)$$

where $\Phi_k$ is the matrix defined in (13), and $y^k$ is the observation given in (10).

The stochastic gradient method [1] is an iterative method in which for a given approximation $u^{k-1}$, the new approximation is obtained by marching in the opposite direction of the gradient of the quadratic form $H_k(u)$ in (42). More precisely, given a step size, $\mu_k > 0$, $u^k$ is computed from

$$u^k = u^{k-1} + \mu_k(b^k - A_k u^{k-1}), \quad \mu_k > 0. \qquad (43)$$

Equation (43) is said to be stochastic gradient method because it is not the classic gradient method for computing (19). A classic gradient method for (19) is an iterative method in which each step is defined by the same matrix $A_k$. In the stochastic gradient method of (43), however, the matrix $A_k$ changes with each step, as shown in Fig 1(b).

The stochastic gradient method (43) is a special case of the canonical form of the constrained approximation of (40) with

$$B = I, \gamma = 0. \qquad (44)$$

We will show in a later section that, under certain conditions, the approximation $u^k$ computed from (43) converges to $\hat{u}$. However, the rate of convergence depends on the condition number $\kappa(A)$ of $A$ in (7); the larger the condition number is, the slower the convergence is. This is evident if $A_k$ is replaced by $A$ in (43), which becomes the well studied classic gradient method. The convergence rate can be improved by the technique of preconditioning.

### B. Preconditioned stochastic gradient method (PSGM)

We now apply preconditioning to the quadratic function in (42). Let $B, C$ be two symmetric $M \times M$ matrices. If $B$ is positive semi-definite, and $C$ is positive definite, then for any $\gamma > 0$, $B + \gamma C$ is positive definite. Define new variable

$$\Psi = (B + \gamma C)^{-1}. \qquad (45)$$

Then minimizing (42) is equivalent to minimizing the following quadratic form

$$H_k^\Psi(w) = (\Phi_k \Psi^{1/2} w - y^k)^T(\Phi_k \Psi^{1/2} w - y^k), w = \Psi^{-1/2} u. \qquad (46)$$

The stochastic gradient method applied to (46) yields

$$w^k = w^{k-1} + \mu_k \Psi^{1/2}(b^k - A_k \Psi^{1/2} w^{k-1}), \quad \mu_k > 0. \qquad (47)$$

Multiplying $\Psi^{1/2}$ from the left, (47) becomes

$$u^k = u^{k-1} + \mu_k (B + \gamma C)^{-1}(b^k - A_k u^{k-1}), \quad \mu_k > 0. \qquad (48)$$

Equation (48) is the preconditioned stochastic gradient method, and the matrix $B + \gamma C$ is said to be the preconditioner. It is now clear that the canonical form of soft constrained approximation (40) is equivalent to the preconditioned gradient method (48), if $B, C$ do not depend on $k$.

It is important to point out that although $A_k$ appears in (40) and (48), it needs not be explicitly computed or stored. All that is needed is a means for computing $A_k u$ for a given vector $u$. On the other hand, the matrices $B, C$ are needed explicitly, but they are normally chosen to be sparse matrices which can be efficiently manipulated. Furthermore, since $B, C$ are independent of $k$, the computation of $(B + \gamma C)^{-1} r^k$ can be implemented efficiently by, e.g., performing an LU factorization of $B + \gamma C$ once for all $k$, and then for each $k$, the forward and backward substitutions are performed on $r^k$.

Next, let $\| \cdot \|_F$ be the Frobenius norm of a matrix. Define

$$d = \|(B + \gamma C)^{-1}\|_F. \qquad (49)$$

### C. Requirement on preconditioning matrix $B$

It is easy to see that $u^k$ from (48) does not converge to $\hat{u}$ even for the seemingly best choices of $\mu, \gamma, B, C$:

$$\mu = 1, \gamma = 0, B = A_k. \qquad (50)$$

Indeed, substitution of (50) into (48) results in (19). If sample sets $x^k$ are drawn independently and identically, and since the solution in (19) only depends on the currently drawn samples, all $u^k$ are expected to have the same statistics, and therefore, it is not expected that $u^k$ would get closer to $\hat{u}$ as $k$ increas-es. Furthermore, it is not known what the expected value $E(u^k) = E(A_k^{-1} b_k)$ is without explicitly knowing the distribution of the random samples $x^k$. This observation provides a heuristic for $B$ to be a constant matrix independent of iteration $k$, which will be assumed in the rest of this paper.

A purpose of using the preconditioner is such that the iteration in (48) has faster convergence. Intuitively, if $A_k$ in (48) is replaced by $A$ of (7), the convergence rate is determined by the condition number of the preconditioned system, and to achieve fast convergence, the condition number needs to be small, i.e.,

$$\kappa((B + \gamma C)^{-1} A) \approx 1. \tag{51}$$

Equation (51) implies, since we can control the parameter of constraint $\gamma$, making it as small as desired, that $B$ needs to be a good approximation of $A$. However, what does it precisely mean by a good approximation, and perhaps more fundamentally, what is the most general condition on $B$ so that PSGM converges?

It turns out that a rather general requirement for $B$ is given by the following

$$(Bx)^T Ax > 0, \text{ for all } x \text{ with } \|x\| = 1. \tag{52}$$

Equation (52) can be interpreted as: if for any unit vector $x$, the angle between the directions $Ax$ and $Bx$ is less than 90 degrees, then $B$ is a good approximation of $A$. Another interpretation of (52) is that matrix $B$ is positive definite relative to matrix $A$. Indeed, requirement (52) is equivalent to $B$ being positive definite in the inner product defined by

$$<x, y>_A = y^T Ax. \tag{53}$$

It can be shown that (52) is equivalent to the following

$$x^T B^{-1} Ax > 0, \text{ for all } x \text{ with } \|x\| = 1. \tag{54}$$

Requirement (54) is met for at least two matrices. For $B = I$, the identity matrix, (54) becomes $x^T Ax > 0$, which is satisfied because $A$ is positive definite. This choice, however, does not provide an improved convergence rate for the stochastic gradient method, but it can be an implementation for imposing a soft constraint represented by $\gamma C$.

For $B = A$, (54) becomes $x^T x > 0$, which is trivially satisfied for unit vector $x$. This would be the ideal choice for the preconditioner if $A$ has favorable properties such as being sparse. Note that using preconditioner $B = A$ in the preconditioned stochastic method of (48) is different from solving the theoretical equation (6). Iteration (48) with $B = A$ can be computed if the probability density function is known, but vector $b$ of (6) cannot be computed as previously stated. Furthermore, the solution $u^k$ from (48) is also different from (19) which is the result of using $B = A_k$ and $\gamma = 0$.

In practice, $A$ may be such that it is costly to compute $(A + \gamma C)^{-1} u$, hence $B = A$ may be a poor choice in terms of complexity. Even if orthogonal basis functions are used, as we discussed in a previous section, $A$ is not a diagonal matrix if the memory effect is included in the approximation. In such a case, the diagonal of $A$ may be used for preconditioning as $B = diag(A)$.

For completeness, we give below an example of symmetric positive definite matrices $A$ and $B$ which do not meet condition (54) (i.e., not every $B$ is a good approximation of $A$):

$$A = \begin{bmatrix} 1 & 2 \\ 2 & 5 \end{bmatrix}, B = \begin{bmatrix} 2 & 1 \\ 1 & 1 \end{bmatrix}, \begin{bmatrix} 1 & 0 \end{bmatrix} B^{-1} A \begin{bmatrix} 1 \\ 0 \end{bmatrix} = -1 < 0. \tag{55}$$

## V. CONVERGENCE ANALYSIS

In this section, we make formal assumptions on random variables, and analyze convergence of PSGM. We will show that under some simple assumptions, PSGM converges to the theoretical best approximation.

### A. Markov process

The approximation $u^k$ in (48) can be considered to be a Markov process defined by stochastic processes $b^k$ and $A_k$, which are in turn defined by random samples drawn from the stochastic process $X(t)$.

Let sample set $x^k$ be drawn at time instances $t_n^k, n = 1, ..., N$. We assume that two different sets of samples are drawn at two different time intervals, i.e.,

$$t_n^k \in [t_1^k, t_N^k], \ n = 1, ..., N, \text{ and } t_N^k < t_1^{k+1}, \text{ for all } k > 0. \tag{56}$$

Then $b^k$ and $A_k$ are formally defined as

$$\begin{aligned} b^k &= [b_i^k]_{M \times 1}, b_i^k = (1/N) \sum_{i=1}^N \phi_i(X(t_n^k)) f(X(t_n^k)), \\ A_k &= [a_{ij}^k]_{M \times M}, a_{ij}^k = (1/N) \sum_{n=1}^N \phi_i(X(t_n^k)) \phi_j(X(t_n^k)). \end{aligned} \tag{57}$$

From definition (57) and (7), it is easy to verify that

$$E(b^k) = b, \ E(A_k) = A. \tag{58}$$

With the definition of random time sampling in (56), we can assume that each of $b^k$ and $A_k$ is an independent and identically distributed (i.i.d.) random process.

Next, we introduce new random variables which are the deviations of $b^k$ and $A_k$ from their means as

$$\theta^k = b^k - b, \ \Omega_k = A_k - A. \tag{59}$$

With $\Psi$ defined by (45), the Markov process (48) can now be written in terms of $\theta^k$ and $\Omega_k$ as

$$u^k = u^{k-1} + \mu_k \Psi(b - Au^{k-1}) + \mu_k \Psi(\theta^k - \Omega_k u^{k-1}). \tag{60}$$

### B. Assumptions

Most of the assumptions on $\theta^k$ and $\Omega_k$ may be derived from those on $b^k$ and $A_k$. Additionally, we assume that the covariances of $\theta^k$ and $\Omega_k$ exist. Write matrix $\Omega_k$ in columns and concatenate all columns into a long vector to define

$$\text{vec}(\Omega_k) = [\omega_1^{kT} \ \cdots \ \omega_M^{kT}]^T, \omega_j^k \text{ is column } j \text{ of } \Omega_k. \tag{61}$$

The covariance of $\Omega_k$ is defined as $E(\text{vec}(\Omega_k) \text{vec}(\Omega_k)^T)$.

We now summarize the assumptions and properties of $\theta^k$ and $\Omega_k$ in the following:

1. $\theta^k$ are i.i.d. with $E(\theta^k) = 0$;
2. $\Omega_k$ are symmetric, i.i.d. with $E(\Omega_k) = 0$;
3. Following covariances exist

$$\sigma_{\theta\theta}^2 \triangleq \| E(\theta^k \theta^{kT}) \|_F,$$
$$\sigma_{\Omega\Omega}^2 \triangleq \| E(\text{vec}(\Omega_k)\text{vec}(\Omega_k)^T) \|_F, \quad (62)$$
$$\sigma_{\Omega\theta}^2 \triangleq \| E(\text{vec}(\Omega_k)\theta^{kT}) \|_F;$$

4. $u^{k-1}$ and its products are pair-wise independent of $\theta^k$ and $\Omega_k$ and their products.

In (62), $\|\cdot\|_F$ denotes the Frobenius norm of a matrix. Equation (62) represents the only assumptions we will need on the random processes $\theta^k$ and $\Omega_k$.

Note that covariances $\sigma_{\theta\theta}^2, \sigma_{\Omega\Omega}^2, \sigma_{\Omega\theta}^2$ in (62) are independent of $k$ because the variables $\theta^k$ and $\Omega_k$ are i.i.d. However, each covariance is a function of the sample size, $N$, as shown in (57). Thanks to the law of large numbers, (57) implies that each covariance can be made arbitrarily small if the number of samples $N$ is large enough.

In (62), assumptions 1-3 are self-explanatory, but assumption 4, regarding the independence of $u^{k-1}$ with the random variables $\theta^k$ and $\Omega_k$, and their products, warrants some discussion. Referring to Fig 1(b), vector $u^{k-1}$ is computed by using samples $x^{k-1}$ drawn in the previous step $k-1$. The outcome of $u^{k-1}$ will not influence in any way how the current set of samples $x^k$ is drawn. The variables $\theta^k$ and $\Omega_k$ are solely dependent on the current sample set $x^k$, the unknown function $f$ and the basis functions $\Phi$, and therefore the outcome of $\theta^k$ and $\Omega_k$, and their products, will not be affected by the outcome of previously computed $u^{k-1}$, which justifies the assumption.

All the results in the rest of this section will be based on the assumptions (54) on $B$, and (62) on $\theta^k$ and $\Omega_k$. Furthermore, we will use $\|\cdot\|_2$ to denote the matrix norm that is induced by the $\ell_2$-norm of vectors [14]. If $\|\cdot\|$ is used for a matrix, it is meant to be the Frobenius norm by default. For vectors, $\ell_2$-norm is always used throughout this paper.

The proofs for all lemmas and theorems of this section are postponed to Appendix.

*C. Mean analysis*

The error in $u^k$ is defined in terms of $\hat{u}$ of (5) as
$$e^k = u^k - \hat{u}. \quad (63)$$

Substituting (60) into (63), it is easy to show that
$$e^k = (I - \mu_k \Psi(A + \Omega_k))e^{k-1} + \mu_k \Psi(\theta^k - \Omega_k \hat{u}). \quad (64)$$

Taking the mean of (64), and noting that $e^{k-1}$ and $\Omega_k$ are independent and that $E(\theta^k) = E(\Omega_k) = 0$, we have
$$E(e^k) = (1 - \mu_k \Psi A)E(e^{k-1}), \quad (65)$$

which leads to following estimate for the norm of the mean

$$\| E(e^k) \| \leq \| 1 - \mu_k \Psi A \|_2 \| E(e^{k-1}) \|. \quad (66)$$

An estimate for the first factor can be found as follows.

**Lemma 1.**

*1) There exists a $\gamma_0 > 0$ such that*
$$\lambda_{\min} \triangleq \min_{0 \leq \gamma \leq \gamma_0} \min_{\|x\|=1} x^T (B + \gamma C)^{-1} A x > 0. \quad (67)$$

*2) Define*
$$\lambda_{\max} \triangleq \max_{0 \leq \gamma \leq \gamma_0} \|(B + \gamma C)^{-1} A\|_2, \tau \triangleq \lambda_{\max} / \lambda_{\min},$$

$$\lambda \triangleq \begin{cases} \lambda_{\max}(1 - \sqrt{1 - \tau^{-2}})\tau, & \text{if } \tau \geq 1 \\ \lambda_{\max}/2, & \text{otherwise} \end{cases}, \quad (68)$$

$$\mu_0 \triangleq 2\frac{\lambda_{\min} - \lambda}{\lambda_{\max}^2 - \lambda^2} = \begin{cases} 1/(\tau\lambda_{\max}), & \text{if } \tau \geq 1 \\ 8(2-\tau)/(3\tau\lambda_{\max}), & \text{otherwise} \end{cases}.$$

*Then the following holds for all $0 \leq \gamma \leq \gamma_0$ and $0 < \mu_k \leq \mu_0$:*
$$\| I - \mu_k(B + \gamma C)^{-1} A \|_2 \leq |1 - \mu_k \lambda|. \quad (69)$$

The key result of Lemma 1 is the existence of $\gamma_0$, $\mu_0$ and $\lambda$ for (69) to hold. The variables defined in (68) are some specific choices to make (69) true, but other definitions are possible too. In particular, it is always possible to choose these parameters such that $0 < \mu_k \lambda \leq 1$. For example, in (68) when $\tau \geq 1$, we have $\mu_k \lambda \leq \mu_0 \lambda = 1 - (1 - \tau^{-2})^{1/2} \leq 1$.

It follows from (66) and (69) that
$$\| E(e^k) \| \leq \| E(e^0) \| \prod_{j=1}^k |1 - \mu_j \lambda|, \quad (70)$$

which leads to the following result.

**Theorem 1.**

*There exist $\gamma_0 > 0$ and $\mu_0 > 0$ such that if the sequence $0 < \mu_k \leq \mu_0$ satisfies $\sum_{k=1}^{\infty} \mu_k = \infty$, then*
$$\lim_{k \to +\infty} E(e^k) = 0, \text{ for all } 0 \leq \gamma \leq \gamma_0. \quad (71)$$

Theorem 1 shows that the mean of preconditioned stochastic gradient iteration converges to the theoretical best approximation. The convergence rate is determined by the factors in (70), i.e., $|1 - \mu_k \lambda|, k = 1, 2, \ldots$. The smaller the factors are the fast the convergence is. Lemma 1 provides some insights into the size of these factors. In most practical problems, e.g., if matrix $(B + \gamma C)^{-1} A$ is symmetric positive definite, we have $\tau \geq 1$ in (68). When $\tau \geq 1$ and if $\mu_k = \mu_0$ in Lemma 1, the factors in (70) are given by
$$|1 - \mu_k \lambda| = 1 - \mu_0 \lambda = 1 - (1 - \sqrt{1 - \tau^{-2}}) = \sqrt{1 - \tau^{-2}}. \quad (72)$$

In the following, we will illustrate the size of (72) by a few examples. In each of the following examples, $\gamma$ is a fixed parameter, and Lemma 1 still holds when $\lambda_{\min}$ and $\lambda_{\max}$ are defined for the fixed $\gamma$, i.e., when $\lambda_{\min}$ and $\lambda_{\max}$ are defined without $\min_{0 \leq \gamma \leq \gamma_0}$ in (67) and $\max_{0 \leq \gamma \leq \gamma_0}$ in (68), respectively, which is what we will assume in the examples. In all examples, $\tau \geq 1$.

First consider the case with $B = A$, $\gamma = 0$. Note that this solution is different from (19) which can be considered as the solution by setting $B = A_k$, $\gamma = 0$. It follows from (67) and (68), $\lambda_{\min} = \lambda_{\max} = \tau = 1, \mu_0 = 1$, and thus, the value in (72) is

0. This implies the mean at each step of the iteration is zero when $\mu_k = \mu_0 = 1$. Therefore, this provides fastest convergence of the mean.

Next, consider $B = I$, $\gamma = 0$. In this case, $\lambda_{\min}$ and $\lambda_{\max}$ are the smallest and the largest eigenvalues of $A$, respectively, and $\tau = \kappa = \text{cond}(A)$, the condition number of $A$. The value of (72) is $(1 - \kappa^{-2})^{1/2}$, which shows that the convergence rate is related to the condition number $\kappa$ of $A$. The larger $\kappa$ is, the slower the convergence is.

Another case is the method of diagonal loading [19] in which $B = A$, $C = I$, $\gamma > 0$. In this case, the value of (72) is

$$\sqrt{1 - \frac{(1 + \gamma \lambda_{\max}^{-1}(A))^2}{(1 + \gamma \lambda_{\min}^{-1}(A))^2}}, \quad \lambda(A) = \text{eigenvalue of } A.$$

If $\gamma$ is small, this is similar to the case of $B = A$, $\gamma = 0$, i.e., the factors of (72) are close to zero. If $\gamma$ is large, this is similar to the case of $B = I$, $\gamma = 0$, i.e., the factors of (72) are close to $(1 - \kappa^{-2})^{1/2}$.

*D. Mean-square analysis*

Next, we consider the variance of the error, which is defined as the autocorrelation of the error $E(e^k e^{kT})$. By using $e^k$ in (64), the variance can be found as follows.

**Lemma 2.**

$$\begin{aligned}
E(e^k e^{kT}) &= (I - \mu_k \Psi A) E(e^{k-1} e^{(k-1)T})(I - \mu_k A\Psi) \\
&+ \mu_k^2 \Psi E((\theta^k - \Omega^k \hat{u})(\theta^k - \Omega^k \hat{u})^T) \Psi \\
&+ \mu_k^2 \Psi E(\Omega_k e^{k-1} e^{(k-1)T} \Omega_k) \Psi \\
&+ \mu_k^2 \Psi E(\Omega_k e^{k-1} \theta^{kT}) \Psi \\
&- \mu_k^2 \Psi E(\Omega_k e^{k-1} \hat{u}^T \Omega_k) \Psi \\
&+ \mu_k^2 \Psi E(\theta^k (e^{k-1})^T \Omega_k) \Psi \\
&- \mu_k^2 \Psi E(\Omega_k \hat{u}(e^{k-1})^T \Omega_k) \Psi.
\end{aligned} \tag{73}$$

We now need to find a bound for the norm of the RHS of (73) in terms of the norm of $E(e^k)$ or $E(e^k e^{kT})$. The first term on the RHS is already an expression of $E(e^k e^{kT})$, and the second term does not involve $e^k$. The remaining five terms need work and they can be estimated using the following result.

**Lemma 3.**

*Let $P, Q$ be two random matrices and $v, w$ be two random vectors. If $v$ and $w$ and their products are independent of $P$ and $Q$ and their products, then*

$$\| E(P v w^T Q) \|_F \leq \| E(v w^T) \|_F \| E(\text{vec}(P)\text{vec}(Q)^T) \|_F. \tag{74}$$

Now each of the last five terms on the RHS of (73) can be written in the form of (74). For example, for the fourth term on the RHS of (73)

$$\Omega_k e^{k-1} \theta^{kT} = \Omega_k e^{k-1} w^T Q, \text{ where}$$
$$w = [1,0,\ldots,0]^T \in \Re^M, Q = [\theta^k, \vec{0},\ldots,\vec{0}]^T \in \Re^{M \times M}. \tag{75}$$

Therefore, since $e^{k-1}$ and $w$ are independent of $\Omega_k$ and $Q$ in (75), Lemma 3 leads to

$$\begin{aligned}
\| E(\Omega_k e^{k-1} \theta^{kT}) \| &= \| E(\Omega_k e^{k-1} w^T Q) \|_F \\
&\leq \| E(e^{k-1} w^T) \|_F \| E(\text{vec}(\Omega_k)\text{vec}(Q)^T) \|_F \\
&\leq \| E(e^{k-1}) \| \| E(\text{vec}(\Omega_k) \theta^{kT}) \|_F \\
&\leq \sigma_{\Omega\theta}^2 \| E(e^{k-1}) \|.
\end{aligned} \tag{76}$$

Similarly, applying (74) to the third term on the RHS of (73), we find

$$\| E(\Omega_k e^{k-1} e^{(k-1)T} \Omega_k) \| \leq \sigma_{\Omega\Omega}^2 \| E(e^{k-1} e^{(k-1)T}) \|_F. \tag{77}$$

In (76) and (77), the variables $\sigma_{\Omega\theta}$ and $\sigma_{\Omega\Omega}$ are defined in (62). For the last term on the RHS of (73)

$$\begin{aligned}
\| E(\Omega_k \hat{u} e^{(k-1)T} \Omega_k) \|_F &\leq \sigma_{\Omega\Omega}^2 \| E(\hat{u} e^{(k-1)T}) \|_F \\
&\leq \sigma_{\Omega\Omega}^2 \| \hat{u} \| \| E(e^{k-1}) \|.
\end{aligned} \tag{78}$$

Working out all the terms in (73), we have the following result.

**Lemma 4.**

*If the condition of Lemma 1 is met, then*

$$\begin{aligned}
\| E(e^k e^{kT}) \| &\leq ((1 - \mu_k \lambda)^2 + \mu_k^2 \sigma_{\Omega\Omega}^2 d^2) \| E(e^{k-1} e^{(k-1)T}) \| \\
&+ \mu_k^2 d^2 (\sigma_{\theta\theta}^2 + 2\sigma_{\Omega\theta}^2 \| \hat{u} \| + \sigma_{\Omega\Omega}^2 \| \hat{u} \|^2) \\
&+ 2\mu_k^2 d^2 (\sigma_{\Omega\Omega}^2 \| \hat{u} \| + \sigma_{\Omega\theta}^2) \| E(e^0) \|.
\end{aligned} \tag{79}$$

In above, $d$ is given by (49). By applying Lemma 4, we can show the following result.

**Lemma 5.**

*There exist $\gamma_0 > 0$, $\lambda_0 > 0$, $\hat{\mu}_0 > 0$, $\delta > 0$ such that if $0 < \mu_k \leq \hat{\mu}_0$ and $0 \leq \gamma \leq \gamma_0$ then for all $k > 0$,*

$$\begin{aligned}
\| E(e^k e^{kT}) \|_F &\leq \| E(e^0 e^{0T}) \|_F \prod_{j=1}^{k}(1 - \mu_j \lambda_0)^2 \\
&+ \delta^2 d^2 \sum_{j=1}^{k} \prod_{i=2}^{j}(1 - \mu_{k-i+2}\lambda_0)^2 \mu_{k-j+1}^2.
\end{aligned} \tag{80}$$

Finally, the convergence of the variance is given below.

**Theorem 2.**

*There exist $\gamma_0 > 0$ and $\hat{\mu}_0 > 0$ such that if the sequence $0 < \mu_k \leq \hat{\mu}_0$ satisfies $\sum_{k=1}^{\infty} \mu_k = \infty$, and $\sum_{k=1}^{\infty} \mu_k^2 < \infty$, then*

$$\lim_{k \to +\infty} E(e^k e^{kT}) = 0, \text{ for all } 0 \leq \gamma \leq \gamma_0. \tag{81}$$

There exists at least one sequence $\mu_k$ for which the condition of Theorem 2 is satisfied. Indeed, the sequence can be defined by

$$\mu_k = \frac{1}{\lambda_0 (k-1) + \hat{\mu}_0^{-1}}, k = 1, 2, \ldots \tag{82}$$

By using the sequence (82), both Theorem 1 an Theorem 2 hold, and therefore, the mean of $u^k$ converges to $\hat{u}$, and the mean of error squared converges to zero. In this sense, we conclude that $u^k$ converges to $\hat{u}$.

The rate of convergence of mean square $E(e^k e^{kT})$ depends on the two terms on the RHS of (80). The first term is same as the convergence rate of the mean. The second term on the RHS of (80) depends on $\mu_k$ and covariances $\sigma_{\theta\theta}^2, \sigma_{\Omega\Omega}^2, \sigma_{\theta\Omega}^2$. Small values for these variables make the iteration fluctuates less. Note that the covariances $\sigma_{\theta\theta}^2, \sigma_{\Omega\Omega}^2, \sigma_{\theta\Omega}^2$ can be made arbitrarily small by using large sample size $N$ in (57). This leads to the conclusion that the larger sample size $N$, the less fluctuation is in the iteration.

## VI. SIMULATION

In this section, we present results from two simulation cases, one for estimating the camera response function (CRF) in image processing and the other for adaptive channel equalization in data communications. These simulations not only demonstrate the effects of using constraints in PSGM but also support the analysis of the previous sections.

### A. Approximation of camera response function

In image processing, camera response functions are nonlinear functions which introduce artifacts such as blurring in images, and in order to remove or reduce such artifacts, it is often necessary to estimate CRF [6][7]. In this simulation, we use a gamma correction to model the CRF and estimate it by a polynomial. The random process is modeled by a Gaussian mixture model which is a common model for pixels of images. More specifically, the CRF is assumed to be given by

$$f(x) = x^\gamma, \gamma = 1/5.5 . \tag{83}$$

Image pixels are normalized to interval $[0,1]$, and they are modeled by two component Gaussian mixture with probability density function (PDF) given by

$$\rho(x) = (\sigma_1 \sqrt{2\pi})^{-1} e^{-(x-\mu_1)^2/\sigma_1^2} + (\sigma_2 \sqrt{2\pi})^{-1} e^{-(x-\mu_2)^2/\sigma_2^2}, \tag{84}$$
$$\mu_1 = .3, \sigma_1 = .01, \mu_2 = .6, \sigma_2 = .007.$$

The basis functions $\Phi(x)$ are chosen to be orthogonal polynomials of degree less than 10, with respect to PDF given in (84). No noise is explicitly introduced in the simulation.

The theoretical best approximation $\hat{u}$ is computed from (6) by using a large number, 50 million, of samples, i.e., about 5 images of 10 megapixels each, and it is shown as the dotted curve in Fig 2. As is expected, the theoretical best approximation derivates from the CRF significantly at the right end of the interval, because there are very few samples near the end of the interval. Constrained approximation can be used to alleviate this problem.

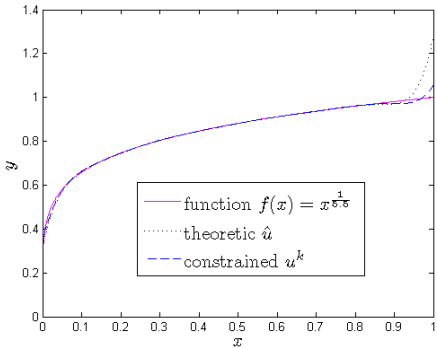

Fig 2 Constrained approximation of CRF. Magenta solid curve is the assumed model for CRF, which is to be approximated. The black dotted curve is the theoretical best approximation as obtained by using 5 million samples. Blue dashed curve is the approximation computed from PSGM.

In the simulation, constraint is implemented by PSGM of (48). The constraint is the requirement that the derivative of the update be smooth, resulting in the following parameters used in the simulation.

$$M = 10, N = 1000, B = I, C = D^T D, \gamma = .02, \mu_0 = .01,$$
$$\mu_k = \mu_0 \text{ if } k < 1000, \mu_k = (k-1000)^{-1}, \text{ if } k > 1000. \tag{85}$$

In (85), $D$ is the difference operator of (31). At each step, 1000 random samples are created from (84), and (48) is used to compute $u^k$. A total of 500,000 steps are performed and the resulting $u^k$ is shown in Fig 2. As is evident, the constrained $u^k$ has smaller deviation from the CRF close to the right end of the interval. This demonstrates the effect of the constrained approximation.

The iteration is performed with a large number of steps to show the convergence history of PSMG, as presented in Fig 3, which show the relative error, $\|u^k - \hat{u}\|/\|\hat{u}\|$ as a function of the iteration number. As is evident, the mean of the error decreases in the iteration, although there are fluctuations in the errors themselves. This is consistent with the analysis of the previous section.

Note that in Fig 2, the computed $u^k$ still differs noticeably from the theoretical $\hat{u}$, especially at the right end of the interval (near 1). This is due to the following two reasons. First, the errors in our analysis are measured by the expected values, $E(e^k)$ and $E(e^k e^{kT})$, both are weighted by the density function $\rho$ of (84), which is very small near 1. Consequently, the seemingly large error in Fig 2 is actually very small when measured by the expected values, as is evident in Fig 3. Secondly, due to round-off errors in simulation, the errors $E(e^k)$ and $E(e^k e^{kT})$ will be bounded below by a nonzero value determined by the machine epsilon of floating point operations, and consequentially $u^k$ may never be the same as $\hat{u}$ no matter how many iterations are performed.

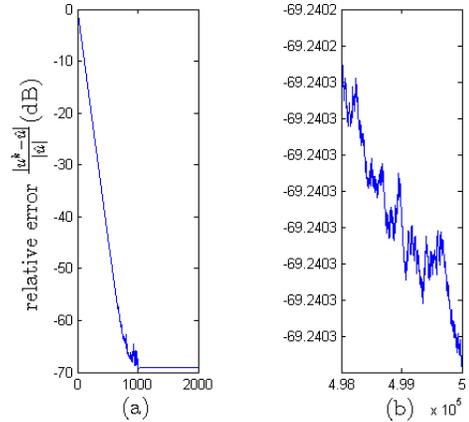

Fig 3 Errors in the iteration. (a) Steps 1 to 2,000. The update factor $\mu_k$ is constant before step 1000, and it is decreasing after step 1000. Although the relative error appears to be constant after step 1000, it is actually decreasing. (b) Steps 498,000 to 500,000 of PSGM

Finally, we point out that error $u^k(x) - f(x)$ may not be a good indication of how well $u^k$ converges to $\hat{u}$ because the error of the theoretical best approximation itself, $\hat{u}(x) - f(x)$, may be large because $f(x)$ may not be well approximated by any combination of the basis functions $\{\Phi_k\}$.

### B. Channel equalization in data communication

In data communication using twisted pair cables, data is transmitted through a channel which introduces distortion into the transmitted data. In this simulation, data from an experimental data transmission system with an RF amplifier is used [20]. The simulation setup is shown in Fig 4.

In Fig 4, function $f$ represents the channel of the twisted pair cable and RF amplifier. The receiver equalizer is to perform an inversion of the channel $f$. In the adaptive equalization, the PSGM is used to approximate the inverse of the channel. Note that while the equalizer (the block named "RX EQ" in Fig 4) itself must be operating in real-time, the adaptation algorithm (the block named "find $u^k$") does not need to be performed real-time, although fast speed is preferred. A small amount, infrequently captured transmitted data, and the corresponding synchronized received data are available to the adaptation algorithm.

In this simulation, since the transmitted and received data are captured experimentally, the function $f$ representing the channel is unknown. Furthermore, the system is set up so that there is a fairly large amount of noise and distortion in the channel.

The channel is nonlinear and has a memory of a very long time duration covering hundreds to a thousand samples. The memory effect may be largely taken care of by a linear filter, and therefore, the strategy used in this simulation is to perform an adaptive linear filter followed by a nonlinear compensation. Although the PSGM is also used for the linear filter, the adaptive linear filter techniques are well studied so it is not of the interest of this paper. In this rest of this paper, we describe the details of nonlinear compensation which follows the linear filter.

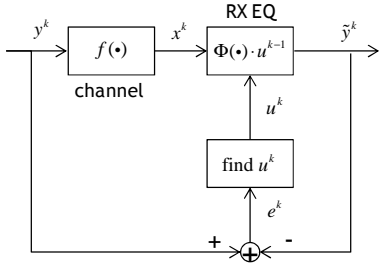

Fig 4 Adaptive channel equalization in data communication. The receiver equalizer is implemented by the PSGM.

*Basis functions*

After the linear filter, the memory effect may be handled by using a small number of taps. In this simulation, the basis functions are based on 5-tap piecewise constant functions. The nonlinear equalizer is implemented by the function given by

$$\Phi(x_n) = \sum_{t=-2}^{2} x_{n-t} \Phi^t(x_{n-t}). \qquad (86)$$

Each of $\Phi^t(x), t=-2,...,2$, is a piecewise constant function represented by a look-up-table (LUT) of 10 bits. That is, there are 1024 constants in each of $\Phi^t(x)$, and each constant is an unknown. Therefore, the total number of unknowns for this simulation is $M = 5120$. These unknowns are the entries of the vector $u$ as discussed in the previous sections, and they will be computed by the PSGM of (48).

*Transmitted and received data samples*

The transmitted and received data are sampled at the rate of 1GHz [20]. The samples are real valued. A large set of samples is captured for transmitted and received data, $x^0$ and $y^0$. Each of them contains $N_0 = 260000$ samples. The transmitted and received samples are synchronized. The samples are scaled to the interval $[-1,1]$.

At each step $k > 0$ of the iteration, a set of transmitted samples $x^k$ and a set of received samples $y^k$ are selected from $x^0$ and $y^0$, respectively. Each set of $x^k$ and $y^k$ contains consecutive time domain samples, but starts from a random location in $x^0$ and $y^0$. The total number of samples in each of $x^k$ and $y^k$ is $N = 1000$.

*Autocorrelation matrix A*

The set of the received samples is used to compute the autocorrelation matrix $A$ of (7). Note that the matrix is not needed in the PSGM process, but it is used for discussion purposes. Although the definition of $A$ requires the probability density function, for the purpose of this simulation, it is approximated by (16) by using the large set of samples, $y^0$. As previously pointed out, even though the basis functions for each of $\Phi^t(x), t=-2,...,2$ are orthogonal, matrix $A$ is not a diagonal matrix because of memory taps. In fact, with the samples captured, $A$ has the sparsity structure, which shows the locations of the nonzero entries of the matrix, as shown in Fig 5.

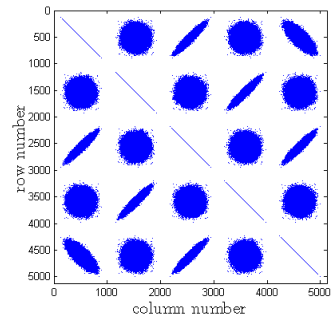

Fig 5 Sparsity of the autocorrelation matrix $A$. Each dot indicates a nonzero entry in the matrix. The matrix is computed from (16) by using 260000 received samples with the 5-tap piecewise constant functions of 10 bits. The dimension of the matrix is $5120 \times 5120$.

*Preconditioning matrix B*

The main diagonal of $A$ is used as the preconditioner in PSGM, i.e., $B = diag(A)$. The use of the preconditioner helps to speed up the convergence by reducing the effective condition number of PSGM. As the condition number of $A$ is very large, Matlab has difficulty to compute it accurately for $M > 80$ (corresponding to 5-tap 4bit LUTs). To demonstrate the effectiveness of the diagonal preconditioner, 4bit LUTs are used to compute the eigenvalues of $A$ and $B^{-1}A$, and the re-

sults are given below:

$$\frac{\lambda_{\max}(A_{80\times 80})}{\lambda_{\min}(A_{80\times 80})} = 5.6\times 10^5, \frac{\lambda_{\max}(B^{-1}A_{80\times 80})}{\lambda_{\min}(B^{-1}A_{80\times 80})} = 7.2\times 10^3. \quad (87)$$

As shown in (87), the preconditioning reduces the condition number of the system by two orders of magnitude.

*Constraint C*

Without using constraints, the computed LUTs $\Phi^t(x)$ may be very jittery. This may be the case even when the number of samples is large. For example, when (19) is solved by using $N_0 = 260000$ samples of $x^0$ and $y^0$, the result is shown in Fig 6 (a). As is evident, the LUTs are discontinuous, and therefore, merely using a large number of samples does not necessarily provide smoothness in the LUTs.

This issue can be resolved by using the smoothness constraint given by matrix $D$, the difference operator of (31). The constraint is incorporated into the preconditioner of the PSGM (48). As a summary, the parameters used in the PSGM simulation is given by

$$\mu_k = 0.1, B + \gamma C = diag(A) + \gamma D^T D, \gamma = 0.02. \quad (88)$$

Note that the preconditioner $B + \gamma C$ of (88) is positive definite and tri-diagonal and hence easy to invert.

The LUTs computed from PSGM (48) are shown in Fig 6 (b). As is evident, the constrained PSGM results in smooth LUTs.

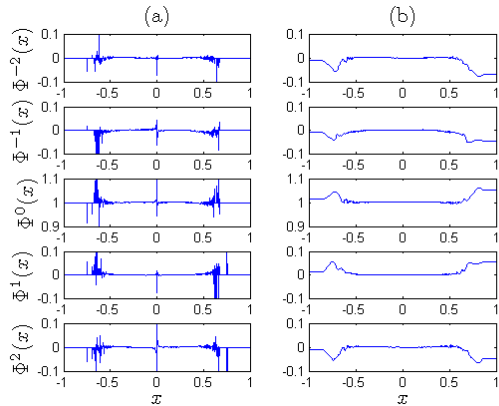

Fig 6 Computed look-up-tables (LUTs). (a) the LUTs are computed from (19) by using 260,000 samples. (b) the LUTs are computed by using constrained PSGM.

*Convergence*

To measure the accuracy of the PSGM, the error is defined as the difference between the transmitted signal and the signal after the receiver equalizer. Referring to Fig 4, at each step $k$ of the iteration, the computed LUTs are used in the RX EQ, with the large set of samples $y^0$ as equalizer input, and the output $\tilde{x}^k$ is compared with large set of the transmitted samples $x^0$. More precisely, the error vector $e^k$ is defined as

$$e^k = [e^k_n]_{N_0\times 1}, e^k_n = x^0_n - \sum_{t=-2}^{2}\Phi^t(x^0_{n-t}), k = 1,.... \quad (89)$$

The size of error as a function of the iteration number is shown in Fig 7. As is consistent with the analysis, the mean of the error decreases, but the errors themselves fluctuate. Note that in this simulation, the update factor $\mu_k = .1$ is a constant through the simulation.

At the last step, $k = 10000$, the computed LUTs are shown in of Fig 6 (b) and the resulting relative error is -41.2dB. As a comparison, when the solution of (19) is computed with all samples of $x^0$ and $y^0$, the LUTs are shown in Fig 6 (a), and the resulting relative error is -41.5dB. Even though the LUTs in Fig 6 (a) is specifically tuned to the data sets $x^0$ and $y^0$, the resulting error is only about 0.3dB better than that from the constrained stochastic process with LUTs of Fig 6 (b), which has much better properties.

Next, the sample plots are shown in Fig 8; they are results before EQ, after linear EQ, and after nonlinear EQ. The after nonlinear EQ result is obtained by using the PSGM. It is evident from Fig 8 that the RX EQ has significantly reduced the channel distortion. It is worth noting that these results were obtained purely based on the samples of the signal waveform, no modulation information is used in the processing. The same technique works for either single carrier QAM signals, or multi-carrier OFDM signals.

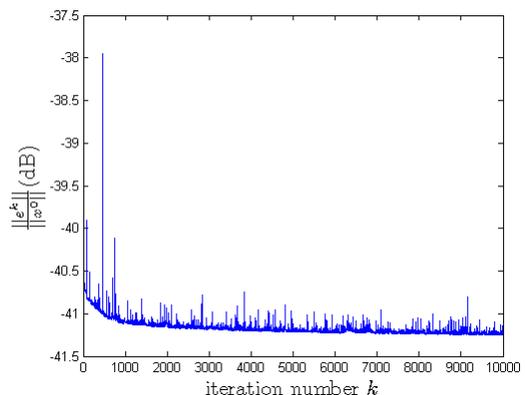

Fig 7 Convergence history of the error computed from PSGM.

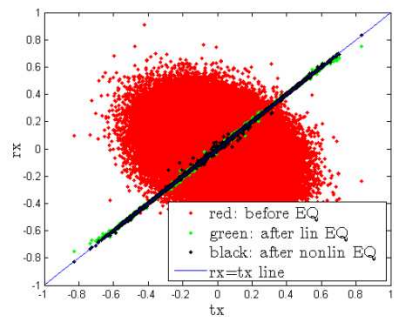

Fig 8 Sample plots. Received samples are plotted against the transmitted samples. Red dots are the sets with received samples before compensated by equalizer, i.e., $x^0$ and $y^0$. Green dots are the sets with received samples after linear EQ. Black dots are the sets with received samples after the linear and nonlinear equalizer, i.e., $x^0$ and $\tilde{y}^k$. The blue line represents the ideal channel in which the received samples equal the transmitted samples. The after nonlinear EQ result is obtained from PSGM.

Finally, we note the reasons why we use a constant update factor $\mu_k = .1$, which meets the requirement of Theorem 1,

but not that of Theorem 2. Data used in this simulation is from a real time prototype system in which the channel is non-stationary. The channel varies with time, often slowly compared to the sample duration, due to such factors as temperature change, clock drifting and component aging etc. Consequently, the stochastic process is not stationary as we assumed at the very beginning of the Section II. However, since the change is slow, it can be regarded as stationary in a short time interval, and the analysis of this paper can be applied within the interval. By using a constant update factor, Theorem 1 implies that the mean of error will approach to zero, so that the channel change can be tracked by the RX EQ. Furthermore, the discussion following Theorem 1 shows that a very small $\mu_k$ leads to slow convergence of the mean, and hence poor tracking of the channel. Therefore, in practical time varying problems, factor $\mu_k$ cannot be set to a too small value since that would lead to the loss of tracking capability. The update factor is often set to a constant because, though the channel varies with time, the signal statistics from one time interval to another are expected to be similar. On the other hand, a constant $\mu_k$ violates Theorem 2, so that the variance of error is no long expected to converge, which is a side effect of the non-stationary stochastic process. This exemplifies the classic trade-off between the speed of tracking (the convergence rate of the mean of the error) and the fluctuation in the result (the variance of the error).

## VII. Conclusion

The constrained and preconditioned stochastic gradient method (PSGM) is investigated. Constraints of the iterative updates can be incorporated into the preconditioner. Both constraining and preconditioning are desirable in stochastic approximations.

Theoretical analysis of the PSGM is performed. Convergence is proved under rather general assumption that the sample sets are drawn identically and independent from one set to another, namely, the assumptions in (62).

Simulations are performed to demonstrate the effect of constraints and to support the results of the theoretical analysis.

## Appendix

**Proof of Lemma 1.**
Let $G(\gamma) = \min_{\|x\|=1} x^T (B + \gamma C)^{-1} Ax$. Then $G(\gamma)$ is a continuous function of $\gamma$. From (54), there exists $\varepsilon' > 0$ so that $G(0) \geq \varepsilon'$. By continuity, there exists $\gamma_0 > 0$ such that for all $0 \leq \gamma \leq \gamma_0$, $G(\gamma) \geq \varepsilon \triangleq \varepsilon'/2$. Therefore,
$$\lambda_{\min} = \min_{0 \leq \gamma \leq \gamma_0} G(\gamma) = \min_{0 \leq \gamma \leq \gamma_0} \min_{\|x\|=1} x^T (B + \gamma C)^{-1} Ax \geq \varepsilon > 0,$$
which proves part 1) of the Lemma 1. Next we show part 2). By definition of $\lambda_{\min}$ and $\lambda_{\max}$, we have $-x^T(B+\gamma C)^{-1}Ax \leq -\lambda_{\min}$ and $x^T A(B+\gamma C)^{-2} Ax \leq \lambda_{\max}^2$. It can be shown that $\lambda$ of (68) satisfies $\lambda < \min\{\lambda_{\min}, \lambda_{\max}\}$. If $\mu_k \leq \mu_0 \triangleq 2 \frac{\lambda_{\min} - \lambda}{\lambda_{\max}^2 - \lambda^2}$, then $-2\lambda_{\min} + \mu_k \lambda_{\max}^2 \leq -2\lambda + \mu_k \lambda^2$.

Therefore, for any $x : \|x\| = 1$, we have
$$\begin{aligned}
&\| (I - \mu_k (B + \gamma C)^{-1} A) x \|_2^2 \\
&= x^T (I - \mu_k (B + \gamma C)^{-1} A)^T (I - \mu_k (B + \gamma C)^{-1} A) x \\
&= 1 - 2\mu_k x^T (B + \gamma C)^{-1} Ax + \mu_k^2 x^T A(B + \gamma C)^{-2} Ax \quad (90) \\
&\leq 1 - 2\mu_k \lambda_{\min} + \mu_k^2 \lambda_{\max}^2 = 1 + \mu_k (-2\lambda_{\min} + \mu_k \lambda_{\max}^2) \\
&\leq 1 + \mu_k (-2\lambda + \mu_k \lambda^2) = (1 - \mu_k \lambda)^2,
\end{aligned}$$
which proves (69) by definition of induced norm of matrix. The value minimizing line 4 of (90) is $\mu_k = \lambda_{\min} / \lambda_{\max}^2$ which is exactly $\mu_0$ of (68) when $\tau \geq 1$.

**Proof of Theorem 1.**
Since $\lambda$ is a positive constant, condition $\sum_{k=1}^{\infty} \mu_k = \infty$ implies $\sum_{k=1}^{\infty} \mu_k \lambda = \infty$, leading to $\prod_{j=1}^{\infty} (1 - \mu_j \lambda) = 0$.

**Proof of Lemma 2.**
The error $e^k$ in (64) can be written as
$$(I - \mu_k \Psi(A + \Omega_k))e^{k-1} - \mu_k \Psi \Omega_k e^{k-1} + \mu_k \Psi(\theta^k - \Omega_k \hat{u}), \quad (91)$$
which is rewritten as, with superscript removed, to save space
$$e^k = \alpha e + \mu \beta e + \mu \eta, \text{ where} \quad (92)$$
$$\alpha = I - \mu_k \Psi(A + \Omega_k), \beta = \Psi \Omega_k, \eta = \Psi(\theta^k - \Omega_k \hat{u}). \quad (93)$$
Then
$$\begin{aligned}
e^k e^{kT} &= (\alpha e + \mu \beta e + \mu \eta)(e^T \alpha^T + \mu e^T \beta^T + \mu \eta^T) \\
&= \alpha e e^T \alpha^T + \mu(\ldots) + \mu^2(\ldots).
\end{aligned} \quad (94)$$
All the terms of (94) linear in $\mu$ can be arranged to have only zero mean variables such as $\Omega_k$ or $\theta^k$, or an additional factor with $e^{k-1}$ which is independent of the zero mean variables. Therefore, after taking the expected value, all the terms linear in $\mu$ disappear and (94) is left with only the first term and the terms that are quadratic in $\mu$, and (94) becomes
$$\begin{aligned}
E(e^k e^{kT}) &= E((I - \mu_k \Psi A)e^{k-1} e^{(k-1)T} (I - \mu_k A \Psi)) \\
&+ \mu_k^2 \Psi E((\theta^k - \Omega_k \hat{u})(\theta^k - \Omega_k \hat{u})^T) \Psi \\
&+ \mu_k^2 \Psi E(\Omega_k e^{k-1} e^{(k-1)T} \Omega_k) \Psi \\
&+ \mu_k^2 \Psi E(e^{k-1}(\theta^k - \Omega_k \hat{u})^T) \Psi \\
&+ \mu_k^2 \Psi E((\theta^k - \Omega_k \hat{u})e^{(k-1)T}) \Psi.
\end{aligned} \quad (95)$$
Equation (95) can be further expanded to (73), which proves Lemma 2.

**Proof of Lemma 3.**
Let $p_j$ be the columns of $P$ and $q_i^T$ the rows of $Q$. Then
$$Pvw^T Q = \sum_{i,j} v_j w_i p_j q_i^T. \quad (96)$$
The $\alpha, \beta$ entry of (96) is given by
$$(Pvw^T Q)_{\alpha\beta} = \sum_{i,j} v_j w_i p_{\alpha j} q_{i\beta}. \quad (97)$$
Since the products of $v$ and $w$ are independent of those of $P$ and $Q$, taking the mean of (97) yields
$$E(Pvw^T Q)_{\alpha\beta} = \sum_{i,j} E(v_j w_i) E(p_{\alpha j} q_{i\beta}). \quad (98)$$

Applying Cauchy-Schwarz inequality to (98) leads to
$$E(Pvw^TQ)^2_{\alpha\beta} \leq \sum_{i,j} E(v_jw_i)^2 \sum_{i,j} E(p_{\alpha j}q_{i\beta})^2$$
$$\leq \|E(vw^T)\|_F^2 \sum_{i,j} E(p_{\alpha j}q_{i\beta})^2. \quad (99)$$

By definition of Frobenius norm, we have
$$\|E(Pvw^TQ)\|_F^2 = \sum_{\alpha,\beta} E(Pvw^TQ)^2_{\alpha\beta}$$
$$\leq \|E(vw^T)\|_F^2 \sum_{\alpha,\beta,i,j} E(p_{\alpha j}q_{i\beta})^2 \quad (100)$$
$$\leq \|E(vw^T)\|_F^2 \|E(\text{vec}(P)\text{vec}(Q)^T)\|_F^2,$$

which concludes the proof.

**Proof of Lemma 4.**
We already have the estimates for the last five terms on the RHS of (73). For the first term, we need to show that for any two $M \times M$ matrices, the following holds
$$\|PQ\|_F \leq \|P\|_2 \|Q\|_F \text{ and } \|PQ\|_F \leq \|P\|_F \|Q\|_2. \quad (101)$$

Indeed, let the columns of $Q$ be $q_j, j=1,...,M$, then
$$\|PQ\|_F^2 = \|P[q_1,...,q_M]\|_F^2 = \sum_j \|Pq_j\|_2^2$$
$$\leq \sum_j \|P\|_2^2 \|q_j\|_2^2 = \|P\|_2^2 \sum_j \|q_j\|_2^2 = \|P\|_2^2 \|Q\|_F^2. \quad (102)$$

From (102), the first term on the RHS of (73) can be found as
$$\|(I - \mu_k \Psi A)\tau(I - \mu_k \Psi A)\|_F \leq \|I - \mu_k \Psi A\|_2^2 \|\Lambda\|_F$$
$$\leq (I - \mu_k \lambda)^2 \|\Lambda\|_F, \quad (103)$$

where $\Lambda = E(e^{k-1}e^{(k-1)T})$. In (103), the result of (69) is used. The second term on the RHS of (73) has the bound
$$\mu_k^2 d^2(\sigma_{\theta\theta}^2 + 2\sigma_{\Omega\theta}^2 \|\hat{u}\| + \sigma_{\Omega\Omega}^2 \|\hat{u}\|^2), \quad (104)$$

where $d$ is defined in (49). Collecting all the terms, we have
$$\|E(e^ke^{kT})\| \leq ((1-\mu_k\lambda)^2 + \mu_k^2\sigma_{\Omega\Omega}^2 d^2)\|E(e^{k-1}e^{(k-1)T})\|_F$$
$$+ \mu_k^2 d^2(\sigma_{\theta\theta}^2 + 2\sigma_{\Omega\theta}^2\|\hat{u}\| + \sigma_{\Omega\Omega}^2\|\hat{u}\|^2) \quad (105)$$
$$+ 2\mu_k^2 d^2(\sigma_{\Omega\Omega}^2\|\hat{u}\| + \sigma_{\Omega\theta}^2)\|E(e^{k-1})\|.$$

Finally, if the condition of Lemma 1 is met, it can be derived from (70) that
$$\|E(e^{k-1})\| \leq \|E(e^0)\|. \quad (106)$$

Substituting (106) into (105), we have proved Lemma 4.

**Proof of Lemma 5.**
It can be easily shown that
$$\text{if } \mu_k \leq \lambda/(3\lambda^2/4 + \sigma_{\Omega\Omega}^2 d^2),$$
$$\text{then } (1-\mu_k\lambda)^2 + \mu_k^2\sigma_{\Omega\Omega}^2 d^2 \leq (1-\mu_k\lambda/2)^2. \quad (107)$$

Therefore, let
$$\hat{\mu}_0 = \min\{4\lambda/(3\lambda^2 + 4\sigma_{\Omega\Omega}^2), 2(\lambda_{\min} - \lambda)/(\lambda_{\max}^2 - \lambda^2)\},$$
$$\delta^2 = \sigma_{\theta\theta}^2 + 2\sigma_{\Omega\theta}^2 \|\hat{u}\| + \sigma_{\Omega\Omega}^2 \|\hat{u}\|^2 \quad (108)$$
$$+ 2(\sigma_{\Omega\Omega}^2\|\hat{u}\| + \sigma_{\Omega\theta}^2)\|E(e^0)\|, \quad \lambda_0 = \lambda/2.$$

We have from (79),
$$\|E(e^ke^{kT})\| \leq ((1-\mu_k\lambda)^2 + \mu_k^2\sigma_{\Omega\Omega}^2 d^2)\|E(e^{k-1}e^{(k-1)T})\|$$
$$+ \mu_k^2 d^2(\sigma_{\theta\theta}^2 + 2\sigma_{\Omega\theta}^2\|\hat{u}\| + \sigma_{\Omega\Omega}^2\|\hat{u}\|^2)$$
$$+ 2\mu_k^2 d^2(\sigma_{\Omega\Omega}^2\|\hat{u}\| + \sigma_{\Omega\theta}^2) \quad (109)$$
$$\leq (1-\mu_k\lambda_0)^2 \|E(e^{k-1}e^{(k-1)T})\| + \mu_k^2\delta^2 d^2,$$

for all $k > 0$, if $\mu_k \leq \hat{\mu}_0$.

Now applying (109) recursively, we have proved Lemma 5.

**Proof of Theorem 2.**
Let $R \triangleq \sum_{k=1}^\infty \mu_k, S \triangleq \sum_{k=1}^\infty \mu_k^2, S_k \triangleq \sum_{i=1}^k \mu_i^2, l_i \triangleq (1-\mu_i\lambda)^2 < 1$, $T_k \triangleq \sum_{j=1}^k \prod_{i=2}^j (1-\mu_{k-i+2}\lambda)^2 \mu_{k-j+1}^2$. We show that if $R = \infty$ and $S < \infty$, then $\lim_{k\to\infty} T_k = 0$. Indeed, for any $K < k$, we have
$$T_k = \mu_k^2 + ... + \prod_{i=K+1}^k l_i\mu_K^2 + \prod_{i=K}^k l_i\mu_{K-1}^2 + ... + \prod_{i=2}^k l_i\mu_1^2$$
$$\leq \mu_k^2 + ... + \mu_K^2 + (\mu_{K-1}^2 + ... + \mu_1^2)\prod_{i=K}^k l_i \quad (110)$$
$$= S_k - S_{K-1} + S_{K-1}\prod_{i=K}^k l_i.$$

Now $\prod_{i=K}^\infty l_i = 0$ for any fixed $K$ because $R = \infty$. Taking limit in (110) leads to $\lim_{k\to\infty} T_k \leq S - S_{K-1}$ for any fixed $K$. Therefore, $\lim_{k\to\infty} T_k \leq \lim_{K\to\infty}(S - S_{K-1}) = 0$.

ACKNOWLEDGMENT

The authors wish to thank Werner Coomans, Keith Chow and Jochen Maes of Bell Labs, Alcatel-Lucent for insightful discussions and for providing experimental data. The authors also thank Bob Wang of Alcatel-Lucent for his interest and discussions. The anonymous reviewers provided insightful comments leading to the improved presentation of this paper, for which we are grateful.

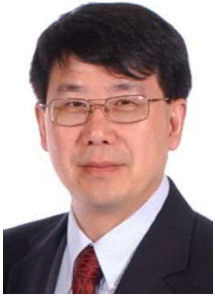

**Hong Jiang** received his B.S. from Southwestern Jiaotong University, M. Math from University of Waterloo, and Ph.D. from the University of Alberta.

Hong Jiang is a researcher and project leader with Alcatel-Lucent Bell Labs, Murray Hill, New Jersey. His research interests include signal processing, digital communications, and image and video compression. He invented key algorithms for VSB demodulation and HDTV video processing in the ATSC system, which won a Technology and Engineering Emmy Award. He pioneered hierarchical modulation for satellite communication that resulted in commercialization of video transmission. He has published more than 100 papers and patents in digital communications and video processing.

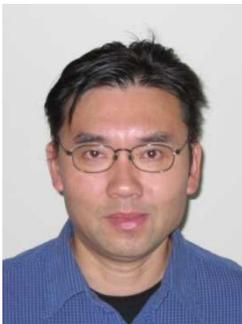

**Gang Huang** received his B.S.E.E. (with highest distinction), and the M.S.E.E. degrees from the University of Iowa, Iowa City, Iowa.

Mr. Huang joined Alcatel-Lucent Bell Labs (formerly AT&T Bell Labs) in 1987. He is currently with the Video Analysis and Coding Research Department of Bell Labs. His current interests include Compressive sensing and its application in image and video acquisition and compression. He holds more than 10 US patents in areas ranging from data networking, data transmission, signal processing, and image processing.

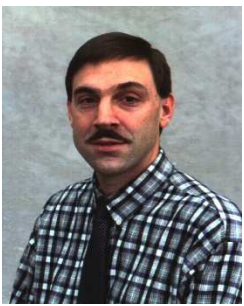

**Paul A. Wilford** received his B.S. and M.S. in Electrical Engineering Cornell University in 1978 and 1979. His research focus was communication theory and predictive coding.

Mr. Wilford is a Bell Labs Fellow, and Senior Director of the Video Analysis and Coding Research Department of Bell Labs. He has made extensive contributions in the development of digital video processing and multi-media transport technology. He was a key leader in the development of Lucent's first HDTV broadcast encoder and decoder. Under his leadership, Bell Laboratories then developed the world's first MPEG2 encoder. He has made fundamental contributions in the high speed optical transmission area.

Currently he is leading a department working on Next Generation video transport systems, hybrid satellite-terrestrial networks, and high-speed mobility networks.

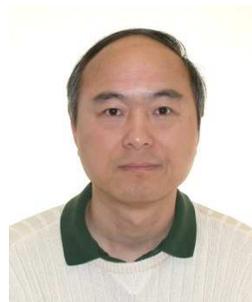

**Liangkai Yu** received his B.S. and M.S. in Electrical Engineering from Fudan University in 1984 and 1987, and Ph.D in Physics from University of Idaho in 1996.

Mr. Yu joined Alcatel Lucent wireless group for 3GPP development in 2010. From 1996, he worked for Lucent Technologies, Agere Systems, Sirius Satellite Radio, and Ibiquity Digital. He was important contributor for the research and development of HDTV VSB demodulator, cable modem, satellite radio receiver, HD Radio, and LTE base station. Currently, his interest is in high speed digital signal processing for wideband radio. His patents and publications were in signal processing and high temperature superconductivity.